\documentclass[11pt]{elsart3-1} 
\usepackage{amssymb}
\usepackage[english,francais]{babel} 
\usepackage[figuresright]{rotating}
\usepackage{amsmath, amsfonts}
\newcommand \RR    {\mathbb{R}}

\newcommand \del   {\partial} 

\newtheorem{theorem}{Theorem}[section]

\newtheorem{e-proposition}[theorem]{Proposition} 

\newtheorem{e-definition}[theorem]{Definition\rm}

\newtheorem{theoreme}{Th\'eor\`eme}[section]

\newtheorem{proposition}[theoreme]{Proposition}

\renewcommand{\theequation}{\arabic{equation}}
\def\og{\leavevmode\raise.3ex\hbox{$\scriptscriptstyle\langle\!\langle$~}}
\def\fg{\leavevmode\raise.3ex\hbox{~$\!\scriptscriptstyle\,\rangle\!\rangle$}}

\newcommand \eps {\epsilon} 
\newcommand \be {\begin{equation}}
\newcommand \ee {\end{equation}}
\newcommand \bei {\begin{itemize}}
\newcommand \eei {\end{itemize}}

\newcommand \la \langle
\newcommand \ra \rangle

\numberwithin{equation}{section}
\newcommand \bel {\be\label}
\newcommand \delb {\overline {\del}} 
\newcommand \Hcal {\mathcal H}
\renewcommand \del \partial
\newcommand \delu {\underline{\del}} 
\newcommand {\delt}{\widetilde{\del}} 
\renewcommand \RR{\mathbb{R}}   
\newcommand {\Mb}{\Mcal_0} 
\newcommand \Mcal {\mathcal M}
\newcommand \Mint {\Mcal^\text{int}}
\newcommand \Mtran {\Mcal^\text{tran}}
\newcommand \Mext {\Mcal^\text{ext}} 

\journal{the Acad\'emie des sciences}
\begin{document} 
\centerline{}
\begin{frontmatter}
 
\selectlanguage{english}
\title{The Euclidian-hyperboidal foliation method 
\\
and the nonlinear stability of Minkowski spacetime}

\selectlanguage{english}
\author[authorlabel1]{Philippe G. LeFloch}
\ead{contact@philippelefloch.org}
\and 
\author[authorlabel2]{Yue Ma}
\ead{yuemath@mail.xjtu.edu.cn}

\address[authorlabel1]{Laboratoire Jacques-Louis Lions \& Centre National de la Recherche Scientifique 
\\
Sorbonne Universit\'e, 4 Place Jussieu, 75252 Paris, France. } 

\address[authorlabel2]{School of Mathematics and Statistics
\\
Xi'an Jiaotong University, Xi'an, 710049 Shaanxi, People's Republic of China.
\\
---
\\
December 2017}

\medskip

\begin{abstract}
\selectlanguage{english}   
We introduce a new method for analyzing nonlinear wave-Klein-Gordon systems and establishing global-in-time existence results for the Cauchy problem 
when the initial data need not have compact support. 
This method, which we call the Euclidian-Hyperboidal Foliation Method (EHFM),  
 relies on the construction of a spacetime foliation obtained by glueing together 
asymptotically Euclidian  
and asymptotically hyperboloidal hypersurfaces. 
Well-chosen frames of vector fields 
(null-semi-hyperboloidal, Euclidian-hyperboloidal) 
 allow us to exhibit the structure of the equations under consideration 
and analyze the decay of solutions in timelike and in spacelike directions. 
New Sobolev inequalities for Euclidian-hyperboloidal foliations
involving the Killing fields of Minkowski spacetime (but not the scaling field), 
as well as pointwise bounds for wave and Klein-Gordon equations on curved spacetimes 
are established. 
Our bootstrap argument involves a hierarchy of (almost sharp) energy and pointwise bounds and distinguishes between low- and high-order derivatives of the solutions. 
We apply this method to the Einstein equations when the matter model 
is a massive  field and the methods by Christodolou and Klainerman and by 
Lindblad and Rodnianski do not apply. 
%
%
%
%
%
%
%
%
%

\end{abstract}
\end{frontmatter}

%
   

\


\selectlanguage{english} 

\tableofcontents


\section{Purpose of this work}

\subsection{Coupled wave-Klein-Gordon systems}

Our main motivation in the present work is the global evolution problem for self-gravitating massive matter fields and, especially, the global nonlinear stability of Minkowski spacetime in this context. 
We introduce here a new method for the global analysis of nonlinear coupled wave-Klein-Gordon systems, which does not require Minkowski's scaling field. 
This method, which we refer to 
as the {\sl Euclidian-Hyperboloidal Foliation Method} (EHFM), 
generalizes the Hyperboloidal Foliation Method (HFM)
 in \cite{PLF-YM-book,PLF-YM-CRAS,PLF-YM-one,PLF-YM-two}. The later method  
extended a methodology first proposed in Klainerman \cite{Klainerman85} and Hormander \cite{Hormander} but also included several novel techniques.
In the present Note, 
we are able to handle solutions that need not be compactly supported, which is essential for the application to the Einstein equations. Our method applies to a broad class of nonlinear wave-Klein-Gordon systems. 

In \cite{PLF-YM-book}, 
we relied on a foliation of the interior of the light cone in Minkowski spacetime by spacelike hyperboloids 
and we  derived sharp pointwise bounds and high-order energy estimates 
in order to study the nonlinear coupling taking place between wave and Klein-Gordon equations.
Our method had the advantage of relying on the Lorentz boosts and the translations only, 
rather than on the full family of (conformal) Killing fields of Minkowski spacetime, as was the case in earlier works such as \cite{CK,LR2}. 

With the Hyperboloidal Foliation Method, we solved the global nonlinear stability problem when the Einstein equations are coupled to a massive scalar field and are expressed in wave gauge; see \cite{PLF-YM-two}. 
However, only the restricted class of initial data sets coinciding with Schwarzschild data outside a spatially compact domain was treated therein. 
  
We provide here an outline of our new method while a full presentation can be found in the Monograph \cite{PLF-YM-c}. 

\subsection{The basic strategy: the interior, transition, and exterior domains, together with adapted frames}

We work in spacetimes $\Mcal$, that is, a 4-manifold endowed with a Lorentzian metric $g$, on which we assume a global foliation $(x^\alpha) = (t,x) = (t, x^a)$ with $\alpha=0,1,2,3$ and $a=1,2,3$. 
Our approach distinguishes within the spacetime 
between interior and exterior spacetime domains, in which different  foliations are required. A gluing technique involving a transition region concentrated near the light cone from the origin is also introduced. 

\bei

\item {\sl Interior domain.} In our approach, 
a region $\Mint \subset \Mcal$ is foliated by spacelike hypersurfaces 
which are constructed from (truncated) hyperboloids in Minkowski spacetime $\RR^{3+1}$. 

\item {\sl Exterior domain.} In our approach, a region $\Mext \subset \Mcal$
is foliated by asymptotically flat hypersurfaces of $\Mcal$, which are constructed from
are flat spacelike hypersurfaces in Minkowski spacetime $\RR^{3+1}$.  

\item {\sl Transition domain.} In a third region $\Mtran \subset \Mcal$, 
he two domains above are glued together around the light cone around which
which the geometry of the foliation changes drastically from begin hyperboloidal to being Euclidian in nature. 

\eei

Several frames of vector fields play a role here: the Cartesian frame $\del_\alpha$, the semi-hyperboloidal frame $\delu_\alpha=(x_a /t)\del_t + \del_a$, as well as the null frame 
$\delt_a =(x^a/r) \del_t+\del_a$  is used
at various stages of our analysis: 

\bei 

\item Vector fields tangent to the foliation: $\delu_a$ in the interior domain
 and $\del_a$ in the exterior domain, 
which are used in expressing the energy estimates. 

\item Vector fields relevant for decomposing the metric and the nonlinearities: $\delu_a$ in the interior domain
 and $\delt_a$ in the exterior domain. 

\eei 
 
\subsection{Main contributions}

A major challenge is to cope with the nonlinear coupling taking place between the geometry and the 
matter terms of various gravity field equations (including the Einstein equations),  
 or equivalently between wave equations and Klein-Gordon equations on a curved spacetime, 
which potentially could lead to a blow-up phenomena of the energy norm under consideration 
and prevent the global existence of solutions.

We rely on basic high-order energy estimates obtained by applying the translations, boosts, and spatial rotations of Minkowski spacetime.

\bei
 
\item Considering a simpler model first and then analyzing the full problem of interest, we provide a {\sl classification of the nonlinearities} arising in 
wave-Klein-Gordon equations 
and systematically we compare them with the terms we control with our energy functional.  

\item By proposing a general and synthetic proof, we establish that our frames of vector fields enjoy {\sl favorable commutator estimates} in order for high-order energy estimates to be derived, 
and, especially, enjoy good commutation properties with the Killing fields of Minkowski spacetime

\item We derive new {\sl Sobolev inequalities for Euclidian-hyperbo\-loidal foliations,} which are established 
by studying cone-like domains first.

\eei 

\bei

\item We establish {\sl sharp sup-norm estimates} for solutions to wave and Klein-Gordon equations on curved spacetime, after decomposing the solution operators in suitable frames and 
building on the following three approaches: an ordinary differential equation argument along rays, a characteristic integration argument, and Kirchhoff's explicit formula.

\eei 

For gravity field problems, the {\sl wave gauge conditions} play a central role in the derivation of 
energy and pointwise bounds 
and especially provide a control on the components of the metric 
associated with a propagation equation that does not satisfy the null condition.
Let us point out that in the case of $1+1$ dimensions discussed in \cite{Ma}, solutions to wave-Klein-Gordon equations 
enjoy much better global properties \cite{Delort01,Stingo}. 
The use of hyperboloidal foliations for wave equations was suggested first by Klainerman \cite{Klainerman85} and Hormander \cite{Hormander}. Earlier on, in \cite{Friedrich81}, Friedrich also studied hyperboloidal foliations of Einstein spacetimes and succeed to establish global existence results for the Cauchy problem for the conformal vacuum field equations.


\section{A new approach: the Euclidian-Hyperboloidal Foliation Method (EFHM)} 

\subsection{A decomposition of the spacetime}

We thus decompose the future $\Mcal$ of an initial hypersurface $\Mcal_0$, 
and 
distinguish between three regions: 
$
\Mcal = \Mint \cup \Mtran \cup \Mext. 
$
Without loss of generality, we label the initial hypersurface as $\big\{ t=1 \big\}$ in our global coordinate chart $(t, x^1, x^2, x^3)$ with $r^2 = \sum (x^a)^2$. We write 
\be
\Mcal \simeq \big\{ t \geq 1 \big\} \subset \RR^{3+1}. 
\ee
Symmetries of  Minkowski spacetime are viewed as approximate symmetries for our spacetime $\Mcal$: 

\bei 

\item The {\sl translations} are generated by the vector fields 
$
\del_\alpha$
$\alpha=0,1,2,3, 
$
which, for instance, will be tangent to the time slices in the exterior domain. 

\item The {\sl Lorentz boosts} are generated by the vector fields 
$
L_a = x_a \del_t  + t \del_a$, 
$a=1,2,3$, 
which, for instance, will be tangent to the time slices in the interior domain. 

\item The {\sl spatial rotations} are generated by the vector fields 
$
\Omega_{ab} = x_a \del_b - x_b \del_a$,
$a=1,2,3$, 
which will be tangent to the time slices in, both, the exterior and the interior domains. 

\eei 

We refer to the above as the family of {\sl admissible vector fieds.}  
These fields commute with the wave and Klein-Gordon operators in Minkowski spacetime.
On the other hand, importantly, throughout our analysis we avoid to rely on {\sl Minkowski's scaling field} 
$
S= t \del_t + r \del_r$, 
since it does not commute with the Klein-Gordon operator in Minkowski spacetime. 

\bei 

\item Within an interior domain $\Mint$, we rely on the foliation based on 
the hyperboloidal slices 
\be
\big\{ t^2 - r^2 =s^2 \big\} \subset \RR^{3+1}
\ee
with hyperbolic radius
 $s \geq 2$. These slices are most convenient in order to analyze wave propagation issues and establish the decay of solutions in timelike directions. 

\item Within an exterior domain denoted by $\Mext$, we rely on Euclidian slices of constant time $c \geq 1$
$$
\big\{ t =c \big\} \subset \RR^{3+1}. 
$$
 These slices are most relevant in order to analyze the asymptotic behavior of solutions  in spacelike directions
and the properties of asymptotically flat spacetimes. 

\eei

\subsection{Time foliation of interest} 
\label{sec:4235}

We take advantage of both foliations above by glueing them together
within a transition region $\Mtran$, as follows.
Consider a cut-off function $\chi=\chi(y)$ (suitably defined to provide a smooth transition between the region) 
satisfying 
$\chi(y) = 0$ for $y \leq 0$, while 
$\chi(y) = 1$ for $y \geq 1$, 
and being increasing within the interval $(0,1)$. Then, we introduce
 the {\sl transition function}  
$$
\xi(s,r) := 1 - \chi(r +1 - s^2/2) \in [0,1],
$$
which is globally smooth and is defined for all $s \geq 1$ and all $r \geq 0$.
The function $\xi$ is not constant precisely in a transition region around the light cone 
$2 r \simeq s^2= t^2 - r^2$.

The {\rm Euclidian-hyperboloidal time function} is the function $T=T(s,r)$ 
defined by the following ordinary differential problem:
\bel{eq 1 foliation}
\aligned
\del_rT(s,r) & = \chi(s-1) \xi(s,r)  \frac{r}{\sqrt{r^2+s^2}},\qquad  
\\
T(s,0) &= s \geq 1.
\endaligned
\ee
By definition, the {\sl Euclidian-Hyperboloidal foliation} is determined from this time function and 
consists of the following family of spacelike hypersurfaces 
$$
\Mcal_s := \big\{ (t,x)\, / \, t = T(s,|x|) \big\}.
$$

Our analyzing is performed in the following spacetime regions: 
\be
\aligned
\Mcal_{[s_0,s_1]} : & = \big\{ (t,x)\, / \, T(s_0, |x|) \leq t \leq T(s_1, |x|)  \big\} \subset \RR^{3+1}, 
\\
\Mcal_{[s_0, +\infty)} :& = \big\{  (t,x)\, / \,  T(s_0,|x|) \leq t \big\}  \subset \RR^{3+1},
\endaligned
\ee
and the interior, transition, and exterior domains (with $r=|x|$) 
\be
\aligned
\Mint_s :& = \big\{ t^2 = s^2+r^2,\quad r\leq -1 + s^2/2 \big\}   &&\text{hyperboloidal region,}
\\
\Mtran_s :& = \big\{  t = T(s,r), \quad - 1 + s^2/2 \leq r\leq s^2/2 \big\}  && \hskip1.cm \text{transition region,}
\\
\Mext_s  :& = \big\{ t= T(s), \quad  r\geq s^2/2 \big\}  &&\text{Euclidian region}.
\endaligned
\ee
By construction, in the interior, the relation  $T^2 = s^2 + r^2$ holds and the slices consist of hyperboloids of Minkowski spacetime. 
On the other hand, 
in the exterior, one has  $T=T(s) \simeq s^2$ which is independent of $r$ and represents a ``slow time'', and the slices consists of flat hyperplanes of Minkowski spacetime.

We choose our cut-off function $\chi$ as follows, by setting first 
\be 
\rho(y) := \left\{
\aligned
&e^{\frac{-2}{1-(2y-1)^2}},\quad &&0 <y< 1, 
\\
&0,\quad &&\text{otherwise}. 
\endaligned
\right.
\ee 
We then set 
$
\rho_0 := \int_{-\infty}^{+\infty}\rho(y) \, dy = \int_0^1\rho(y) \, dy >0
$
and define
\be
\chi(y) := \rho_0^{-1}\int_{-\infty}^y \rho(y') \, dy', \quad y \in \RR. 
\ee 


\subsection{Weighted energy norm}

We use the energy norm associated to the wave and Klein-Gordon equations and induced on 
our Euclidian-hyperboloidal slices. In addition, 
in the exterior domain we introduce a weight function which provides us with the required control of the decay in spacelike directions. This weight depends upon the variable 
$q := r-t$, 
that is, the distance to the light cone from the origin and so, for some $\eta \in (0,1]$ and using our cut-off function, we set 
\be
\omega_\eta (t,r) := \chi(q) (1+q)^{\eta}  = \chi(r-t) ( 1+ r-t)^{\eta}. 
\ee
To any function $v$ defined in the domain $\Mcal_{[s_0,s_1]}$ limited by two slices of our foliation (say $1 \leq s_0 \leq s \leq s_1$), 
we define the {\sl energy functional}  
$$
\aligned
& E_{\eta,c}(s,v)
\\
& := 
\int_{\Mcal_s} (1+\omega_\eta)^2
\Big( 
\Big( 1 -  \chi(1-s)^2 \xi^2 {r^2 \over t^2} \Big) \big( \del_t v \big)^2 
+  \sum_a \Big( \xi {x^a \over t} \del_t v + \del_a v \Big)^2 + c^2 v^2
\Big) dx, 
\endaligned
$$
in which, by definition, one has $t=T(s,r)$ and $r=|x|$ on $\Mcal_s$. 
With the notation  
\be
\aligned
\delb_a 
& := \del_aT(s,r) \del_t + \del_a  
\\
& =  \chi(1-s)  \xi(s,r) {x^a \over T(s,r)} \del_t v + \del_a, 
\endaligned
\ee
we obtain the alternative form 
\be
\aligned
E_{\eta,c}(s,v) 
& = \int_{\Mcal_s}(1+\omega_\eta)^2\Big(|\zeta(s,r)\del_t v|^2 + \sum_a|\delb_a v|^2 + c^2 v^2\Big) \, dx, 
\endaligned
\ee
in which the coefficient $\zeta= \zeta(s,r) \in [0,1]$ is defined as 
\be
\aligned
\zeta & := \chi(1-s) \sqrt{\frac{s^2+\chi^2(r - 1+s^2/2 )r^2}{s^2+r^2}}.
\endaligned
\ee
This energy (together with its generalization to a curved metric) leads us to a control of the  
weighted wave-Klein-Gordon energy associated with the operator $\Box v - c^2 v$ with $c \geq 0$.


\subsection{A Sobolev inequality for the transition and exterior domains}

\begin{proposition} 
\label{prop 2 04-11-2017}
For all sufficiently regular functions defined on $\Mcal_{[s_0,s_1]}$ with $2 \leq s_0 \leq s \leq s_1$, one has: 
\bel{ineq 2 sobolev}
|u(x)|\leq C(1+r)^{-1}\sum_{|I|+|J|\leq 2}\|\delb_x^I\Omega^J\|_{L^2(\Mtran_s\cup \Mext_s)}, 
\qquad x\in \Mtran_s \cup \Mext_s, 
\ee
\bel{ineq 1 sobolev}
|u(x)|\leq C(1+r)^{-1}\sum_{|I|+|J|\leq 2}\|\del_x^I\Omega^J\|_{L^2( \Mext_s)}, 
\qquad x\in \Mext_s. 
\ee
Here, $\delb_x^I$ denotes any $|I|$-order operator determined from the fields $\{\delb_a\}_{a=1,2,3}$, 
while $\del_x^I$ denotes 
any a $|I|$-order operator determined from the fields $\{\del_a\}_{a=1,2,3}$.
\end{proposition}

\

\noindent 
{\bf Proof.} We only sketch the argument. 
We consider the parametrization $(s,r)$ of $\Mcal_{[s_0,+\infty)}$, and 
recall that on $\Mtran_s \cup \Mext_s$, $s$ is constant and $t = T(s,r)$, $-1+s^2/2 \leq r$. We consider the restriction of $u$ on $\Mtran_s \cup \Mext_s$, that is, the function 
$
v_s(x) = u(T(s,r),x)
$
and we remark the relations
$$ 
\aligned
\del_a v_s 
& = \delb_au = (x^a/r){\del t \over \del r}  \del_t u +  \del_a u = \frac{\xi_s(r)x^a}{t}\del_t u+\del_a u,
\\ 
\del_b\del_a v_s 
& = \delb_b\delb_a u,
\\
\Omega_{ab}v_s 
& = (x^a\del_b-x^b\del_a)u = \Omega_{ab}u.
\endaligned
$$
Now we apply a Sobolev inequality for $v_s$ adapted to a cone-like region (cf.~\cite{PLF-YM-c}  for details), 
and  \eqref{ineq 2 sobolev} is established, 
while 
\eqref{ineq 1 sobolev} is established in the same manner.
\qed


\subsection{A Sobolev inequality for the interior domain}

\begin{proposition} 
For all sufficienty regular functions defined in a neighborhood of the hypersurface $\Mint_s$, 
the following estimate holds  
\bel{ineq 3 sobolev}
t^{3/2} \, |u(x)|\lesssim
\sum_{|J| \leq 2} \| L^J u\|_{L^2(\Mint_s)}, 
\qquad x\in\Mint_s.
\ee 
\end{proposition}

\

\noindent {\bf Proof.} We only skecth the proof and refer to~\cite{PLF-YM-c} for details. 
We consider the restriction $v_s(x) := u(\sqrt{s^2+r^2},x)$
of the function $u$ on the hyperboloid $\Hcal_s$ with $|x|\leq - 1 + s^2/2 $. 
Then we have 
$$
\del_av_s = \delb_a u = t^{-1}L_a u = (s^2+r^2)^{-1/2}L_a u.
$$
Take a $x_0\in\Hcal_s$, with out loss of generality, we can suppose that $x_0 = -3^{-1/2}(r_0,r_0,r_0)$. We consider the positive cone $C_{s/2,x_0}\subset \{|x|\leq -1+s^2/2 \}$.
In this cone we introduce the change of variable:
$
y^a: = s^{-1}(x^a-x_0^a)
$
and we define
$
w_{s,x_0}(y): = v_s(sy + x_0)$
for $y \in C_{1/2,0}$. 
Therefore, we obtain 
$$
\aligned
\del_aw_{s,x_0} & = s\del_av_s = \frac{s}{\sqrt{s^2+r^2}}L_au,
\\
\del_b\del_aw_{s,x_0} & = \frac{s^2}{s^2+r^2}L_bL_au -s^2x^b(s^2+r^2)^{-3/2}L_au.
\endaligned
$$
Thus,  for $|I|\leq 2$, we obtain 
$
|\del^I w_{s,x_0}|\leq C\sum_{|J|\leq 2}|L^J u|.
$
Then by a Sobolev inequality adapted to cone-like region (cf~\cite{PLF-YM-c})
 we see that
$$
\aligned
|w_{s,x_0}(0)|^2 \leq C\sum_{|I|\leq 2}\int_{C_{1/2,0}}|\del^Iw_{s,x_0}|^2dy = Cs^{-3}\sum_{|J|\leq 2}\int_{C_{1/2,0}}|L^J u|^2dx,
\endaligned
$$
which leads to
\bel{eq 1 02-11-2017}
|u(x_0)|\leq Cs^{-3/2}\|L^J u\|_{L^2(\Hcal_s)}.
\ee

On the other hand, when $r_0\geq 1 $, we consider the cone $C_{r_0/2,x_0}$ and we introduce the function
$$
w_{x_0}(y):=v_s(r_0y+x_0),\quad y\in C_{1/2,0}.
$$
It is clear that
$
\del_aw_{x_0} = \frac{r_0}{\sqrt{r^2+s^2}}L_au
$
and
$$
\del_b\del_aw_{x_0} = \frac{r_0^2}{r^2+s^2}L_bL_au - r_0^2x^b(s^2+r^2)^{-3/2}L_au.
$$

In the cube $C_{r_0/2,x_0}$ one has $r\geq \frac{\sqrt{3}}{2}r_0$. Thus for $|I|\leq 2$ we find 
$$
|\del^I w_{x_0}|\leq \sum_{|J|\leq 2}|L^Ju|.
$$
Then, by a Sobolev inequality adapted to a cone-like region \cite{PLF-YM-c}, 
 we have 
$$
\aligned
u(x_0)|^2 = |w_{x_0}(0)|^2 
& \leq C\sum_{|J|\leq 2}\int_{C_{1/2,0}}|\del^Iw_{x_0}|^2dy 
\\
& = Cr_0^{-3}\sum_{|J|\leq 2}\int_{C_{1/2,0}}|L^J u|^2dx, 
\endaligned
$$
which leads us to
$
|u(x_0)|\leq Cr_0^{-3/2}\sum_{|J|\leq 2}\|L^J u\|_{L^2(\Hcal_s)}.
$
When $r_0\leq 1$, we have $\sqrt{s^2+r_0^2}\leq 2s$ and the desired result is proved. 
\qed

 \section{Existence theory and application to the Einstein equations}

The proposed method leads us to new results of global-in-time existence for coupled systems of nonlinear wave-Klein-Gordon equations. A particularly interesting model is provided by the coupling of a massive scalar field to the Einstein equations, that is, 
\bel{MainPDE-limit}
\aligned
\Box_g g_{\alpha\beta}
= \, &  F_{\alpha\beta}(g, \del g)
    + 16 \pi \, \big( - \del_\alpha \phi \del_\beta \phi + U(\phi) \, g_{\alpha \beta} \big),
\\
\Box_g \phi  - U'(\phi) = \, &  0, 
\endaligned
\ee
whose unknowns are a Lorentz metric $g=(g_{\alpha\beta})$ (assumed to satisfy the wave gauge conditions)
and a scalar field $\phi$, while 
$U(\phi) = c^2 \phi^2/2$. 
We refer to \cite{PLF-YM-c} for a precise statement of our asymptotic conditions in spatial directions
which allow us to define the notion of $(\sigma, \eps,N)$--asymptotically tame initial data set. Here the exponent $\sigma \in (1/2,1]$ measures the pointwise decay $r^{-\delta}$ of the metric in spacelike directions, in comparison to the Minkowski metric denoted by $g_M$. The following result is a generalization to self-gravitating massive matter of the noninear stability theory in \cite{CK,LR2}.  

\

\begin{theorem}[Nonlinear stability of Minkowski space for self-gravitating massive fields]
\label{theo:main1} 
For all sufficiently small $\eps>0$ and all sufficiently large integer $N$, one can find 
$\sigma, \eta >0$ depending upon $(\eps, N)$  
so that the following property holds. 
Consider any 
$(\sigma, \eps,N)$--asymptotically tame initial data set  
for the Einstein-massive field system \eqref{MainPDE-limit}, i.e. 
$$
(\Mb \simeq \RR^3,g_0, k_0, \phi_0, \phi_1), 
$$
which is assumed to be sufficiently close to a flat and vacuum spacelike slice of Minkowski spacetime $(\RR^{3+1}, g_M)$
in the sense that, in a global coordinate chart $x=(x^a)$ with $r=|x|$, 
\bel{eq:F94}
\aligned 
\|(1+r)^{\eta+|I|} Z^I \del (g_0-g_{M,0}) \|_{L^2(\RR^3)}  &\leq \eps, 
\\
\|  (1+r)^{\eta+|I|} Z^I k_0 \|_{L^2(\RR^3)}  &\leq \eps,
\\
\|(1+r)^{\eta+|I|+1/2} Z^I  \phi_0  \|_{L^2(\RR^3)}   &\leq \eps,
\\
\|(1+r)^{\eta+|I|+1/2} Z^I \phi_1  \|_{L^2(\RR^3)}  &\leq \eps,
\endaligned
\ee 
for all $Z^I =  \del^{I_1} \Omega^{I_2}$ with $|I_1| + |I_2|\leq N$, 
in which $g_{M,0} =  (\delta_{ab})$. 
Then, the corresponding initial value problem associated with \eqref{MainPDE-limit}
admits a global solution
and therefore determines
a globally hyperbolic Cauchy development $(\Mcal, g, \phi)$,
which 
is future causally geodesically complete
 and asymptotically approaches Minkowski metric.
\end{theorem}
 
We recall that a future causally geodesically complete spacetime, by definition, has the property that every affinely parameterized geodesic (of null or timelike type) can be extended toward the future (for all values of its affine parameter).
In view of the theorem for compactly supported kinetic fields established in \cite{FJS3}, 
it would be very natural to also extend our new method to the Einstein-Vlasov model.


\end{document}